\documentclass{amsart}
\usepackage{latexsym}
\usepackage{amssymb}
\usepackage{amsmath}
\usepackage{amsthm}
\usepackage{hyperref}
\usepackage{xcolor}

\newcommand{\CC}{{\mathbb C}}

\newcommand{\RR}{{\mathbb R}}

\newcommand{\SSS}{{\mathbb S}}

\newcommand{\grad}{\mathrm{grad}}

\newcommand{\scal}{\mathrm{Scal}}

\newcommand{\tr}{\mathrm{tr}}

\newcommand{\lgra}{\longrightarrow}

\newtheorem{example}{Examples}[section]
\newtheorem{thm}{Theorem}[section]
\newtheorem{lemma}[thm]{Lemma}
\newtheorem{prop}[thm]{Proposition}
\newtheorem{cor}[thm]{Corollary}
\newtheorem{remark}[thm]{Remark}
\newtheorem{remarks}[thm]{Remarks}
\newtheorem{definition}[thm]{Definition}
\newtheorem{notation}[thm]{Notation}
\newtheorem{exabout:ample}[thm]{Example}

\begin{document}
\title{Biharmonic submanifolds of generalized space forms}
\author[J. ROTH]{Julien Roth}
 \author[A. UPADHYAY]{Abhitosh Upadhyay}
\address[J. ROTH]{Laboratoire d'Analyse et de Math\'ematiques Appliqu\'ees, UPEM-UPEC, CNRS, F-77454 Marne-la-Vall\'ee}
\email{julien.roth@u-pem.fr}
\address[A. UPADHYAY]{Harish Chandra Research Institute, Chhatnag Road, Jhunsi, Allahabad, India, 211019}
\email{abhi.basti.ipu@gmail.com, abhitoshupadhyay@hri.res.in}
\maketitle

\begin{abstract}
We consider biharmonic submanifolds in both generalized complex and Sasakian space forms. After giving the biharmonicity conditions for submanifolds in these spaces, we study different particular cases for which we obtain curvature estimates. We consider curves, complex and Lagrangian surfaces and hypersurfaces for the generalized complex space form as well as hypersurfaces, invariant and anti-invariant submanifolds in case of generalized Sasakian space form. 
\end{abstract}

\section{\textbf{Introduction}}
A harmonic map $\psi$ between two Riemannian manifolds $(M,g)$ and $(N,h)$ is defined as a critical point of the energy functional
$$E(\psi)=\frac{1}{2}\int_M|d\psi|^2dv_g.$$
\indent
In \cite{ES}, Eells and Sampson gave a natural generalization of harmonic maps. A map $\psi$ is called {\it biharmonic} if it is a critical point of the bi-energy functional
$$E_2(\psi)=\frac{1}{2}\int_M|\tau(\psi)|^2dv_g,$$ where $\tau(\psi)$ is the tension field, which vanishes precisely for harmonic maps. G.Y. Jiang \cite{Ji}, studied the first and second variation formulas
of $E_{2}$, which critical maps are called biharmonic maps. The Euler-Lagrange equation associated with this bi-energy functional is $\tau_2(\psi)=0$, where $\tau_2(\psi)$ is the so-called bi-tension field given by 
$$\tau_2(\psi)=\Delta\tau(\psi)-\tr\big(R^N(d\psi,\tau(\psi))d\psi\big).$$
Here, $\Delta$ is the rough Laplacian acting on the sections of $\psi^{-1}(TN)$ and $R^N$ is the curvature tensor of $N$. We will use the following sign convention, i.e.,
$$\Delta V=\tr(\nabla^2V)\quad\text{and}\quad R^N(X,Y)=[\nabla^N_X,\nabla^N_Y]-\nabla^N_{[X,Y]}$$
for any $V\in\Gamma(\psi^{-1}(TN))$ and $X,Y\in\Gamma(TN)$. It is obvious that any harmonic map is biharmonic. There have been extensive studies in this area (see \cite{CMO1, JI, SH, MC, YO, YO1, HU}, for instance). We will focus here on  biharmonic maps which are not harmonic. They are called {\it proper biharmonic maps}.\\ \\
\indent
If the map $\psi: (M, g)\rightarrow (N, h)$ is an isometric immersion from a manifold  $(M,g)$ into an ambient manifold $(N,h)$ then $M$ is called {\it biharmonic submanifold} of $N$. In other words, {\it proper biharmonic submanifolds} are the biharmonic submanifolds which are not harmonic.\\ \\
\indent
In the last decades, biharmonic submanifolds has become a popular subject of research with many significant progresses made by geometers around the world. One of the fundamental problems in the study of biharmonic submanifolds is to know its geometry in space forms. So far, most of the work done has been focused on study of biharmonic submanifolds of space forms.  Several results have been proven in different ambient spaces like real space forms \cite{CMO},  complex space forms \cite{FLMO}, 3-dimensional homogeneous manifolds \cite{LM}, Sasakian space forms \cite{FO}, or products of real space forms (\cite{FOR, Rot}). \\ \\
One of the main problem in the study of biharmonic submanifolds is the Chen's conjecture  \cite{Ch2}:\linebreak
\begin{center}
{\it ``\textbf{The only biharmonic submanifolds of Euclidean spaces are the minimal ones.}''}
\end{center}
This conjecture has been proven in many particular cases (see \cite{Ch3} and references therein for an overview) but it is still open in general. On the other hand, the generalized Chen's conjecture replacing Euclidean spaces by Riemannian manifolds of non-positive sectional curvature turns out to be false (see \cite{LO,OL} for counter-examples). Nevertheless, this generalized conjecture is true in various situations and non-existence results in non-positive sectional curvature is still an interesting question. We will give in this paper two new contexts where such results hold.\\ \\
The purpose of the present paper is to consider the case of generalized space forms. First, we study the case of generalized complex space forms (Section \ref{sec3}). These spaces are Einstein non-K\"ahler Hermitian manifolds which are generalizations of complex space forms. They appear only in dimension 4. After giving the biharmonicity condition for submanifolds of these spaces (see Theorem \ref{thm1}), we consider different particular cases, namely curves, Lagrangian or complex surfaces and hypersufaces. Then, we obtain curvatures estimates as well as non-existence results for biharmonic submanifolds for generalized complex spaces forms of negative (constant) scalar curvature.\\ \\
\indent
In Section \ref{sec4}, we consider generalized Sasakian space forms. This familiy of almost contact metric manifolds generalizes Sasakian space forms and contains also the so-called Kenmotsu and cosymplectic space forms. In this case, we give the general conditions for biharmonic submanifolds (Theorem \ref{thm2}) with a focus on many particular cases such as hypersurfaces and invariant or anti-invariant submanifolds which are the analogous in the contact setting  of complex or totally real (in particular Lagrangian) submanifolds appearing in the complex case. Further, we also obtain some curvature estimates for biharmonic submanifolds of generalized Sasakian space form and as corollary, some non-existence results in case of negative (or appropriately bounded by above) $\phi$-sectional curvature.
\section{\textbf{Preliminaries}}\label{sec2}
\subsection{Generalized complex space forms and their submanifolds}

We begin by giving some basic information about generalized complex space forms. They form a particular class of \linebreak Hermitian manifolds which has not been intensively studied. In 1981, Tricelli and Vanhecke \cite{TV} introduced the following generalization of the complex space forms ($\CC^n$, $\CC P^n$ and $\CC H^n$).\\
\indent
Let $(N^{2n},g,J)$ be an almost Hermitian manifold. We denote the generalized curvature tensors by $R_1$ and $R_2$ which is defined as
$$R_1(X,Y)Z=g(Y,Z)X-g(X,Z)Y,$$
$$R_2(X,Y)Z=g(JY,Z)JX-g(JX,Z)JY+2g(JY,X)JZ, \hspace{.2cm}\forall\hspace{.2cm} X,Y,Z\in \Gamma(TN).$$
The manifold $(N,g,J)$ is called {\it generalized complex space form} if its curvature tensor $R$ has the following form
$$R=\alpha R_1+\beta R_2,$$
where $\alpha$ and $\beta$ are smooth functions on $N$. The terminology comes obviously from the fact that complex space forms satisfies this property with constants $\alpha=\beta$ . \\
\indent
In the same paper \cite{TV}, Tricelli and Vanhecke showed that if $N$ is of (real) dimension $2n\geq6$, then $(N,g,J)$ is a complex space form. They also showed that $\alpha+\beta$ is necessarily constant. This implies that $\alpha=\beta$ are constants in dimension $2n\geq6$, but this in not the case in dimension $4$. Hence, the notion of generalized complex space form is of interest only in dimension $4$.  Further, Olszak \cite {Ols} constructed examples in dimension $4$ with $\alpha$ and $\beta$ non-constant. These examples are obtain by conformal deformation of B\"ochner flat K\"ahlerian manifolds of non constant scalar curvature.  Examples of B\"ochner flat K\"ahlerian manifolds can be found in \cite{Der}. From now on, we will denote by $N(\alpha,\beta)$ a (4-dimensional) generalized complex space form with curvature given by $R=\alpha R_1+\beta R_2$. Note that these spaces are Einstein, with constant scalar curvature equal to $12(\alpha+\beta)$. Of course, they are not K\"ahlerian because if they were, they would be complex space forms. \\ \\
\indent
Now, let $M$ be a submanifold of the generalized complex space form $N(\alpha,\beta)$. The almost complex structure $J$ on $N(\alpha,\beta)$ induces the existence of four operators on $M$, namely
$$j:TM\lgra TM,\ k:TM\lgra NM,\l:NM\lgra TM\ \text{and}\ m:NM\lgra NM ,$$
defined for all $X\in TM$ and all $\xi\in NM$ by
\begin{eqnarray}\label{relationfhst}
JX=jX+kX\quad\text{and}\quad
J\xi=l\xi+m\xi.
\end{eqnarray} 
Since $J$ is an almost complex structure, it satisfies $J^2=-Id$ and for $X,Y$ tangent to $N(\alpha,\beta)$, we have $g(JX,Y)=-g(X,JY)$. Then, we deduce that the operators $j,k,l,m$ satisfy the following relations
\begin{align}
&j^2X+lkX=-X,& \label{relation1.1}\\
&m^2\xi+kl\xi=-\xi,& \label{relation1.2}\\
\label{relation1.3}
&jl\xi+lm\xi=0,&\\
\label{relation1.4}
&kjX+mkX=0,&\\
\label{relation1.5}
&g(kX,\xi)=-g(X,l\xi),&
\end{align}
for all $X\in\Gamma(TM)$ and all $\xi\in\Gamma(NM)$. Moreover $j$ and $m$ are skew-symmetric.

\subsection{Generalized Sasakian space forms and their submanifolds}
Now, we give some recalls about almost contact metric manifolds and generalized Sasakian space forms. For more details, one can refer to (\cite{ABC,Bla,YK}) for instance. A Riemannian manifold $\widetilde{M}$ of odd dimension is said almost contact if there exists globally over $\widetilde{M}$, a  vector field $\xi$, a $1$-form $\eta$ and a field of $(1,1)$-tensor $\phi$ satisfying the following conditions:
\begin{equation} \eta(\xi)=1\quad\text{and}\quad \phi^2=-Id+\eta\otimes\xi.
\end{equation}
Remark that this implies $\phi\xi=0$ and $\eta\circ\phi=0$. The manifold $\widetilde{M}$ can be endowed with a Riemannian metric $\widetilde{g}$ satisfying
\begin{equation}
\widetilde{g}(\phi X,\phi Y)=\widetilde{g}(X,Y)-\eta(X)\eta(Y)\quad\text{and}\quad \eta(X)=\widetilde{g}(X,\xi),
\end{equation}
for any vector fields $X,Y$ tangent to $\widetilde{M}$. Then, we say that $(\widetilde{M},\widetilde{g},\xi,\eta,\phi)$ is an almost contact metric manifold. Three class of this family are of particular interest, namely, the Sasakian, Kenmotsu and cosymplectic manifolds. We will give some recalls about them. \\ \\
\indent
First, we introduce the fundamental $2$-form (also called Sasaki $2$-form) $\Omega$ defined for $X,Y\in\Gamma(TM)$ by 
$$\Omega(X,Y)=\widetilde{g}(X,\phi Y).$$
We consider also $N_{\phi}$, the Nijenhuis tensor defined by
$$N_{\phi}(X,Y)=[\phi X,\phi Y]-\phi[\phi X,Y]-\phi[X,\phi Y]+\phi^2[X,Y],$$
for any vector fields $X,Y$. An almost contact metric manifold is said normal if and only if the Nijenhuis tensor $N_{\phi}$ satisfies
$$N_{\phi}+2d\eta\otimes\xi=0.$$ 
An almost contact metric manifold is said {\it Sasakian manifold} if and only if it is normal and  $d\eta=\Omega$. This is equivalent to
\begin{equation}
(\nabla_X\phi)Y=\widetilde{g}(X,Y)\xi-\eta(Y)X, \hspace{.2 cm}\forall \hspace{.2 cm}X, Y\in\Gamma(\widetilde{M}).
\end{equation}
It also implies that 
\begin{equation}
\nabla_X\xi=-\phi(X).
\end{equation}
An almost contact metric manifold is said {\it Kenmotsu manifold} if and only if $d\eta=0$ and $d\Omega=2\eta\wedge\Omega$. Equivalently, this means
\begin{equation}
(\nabla_X\phi)Y=-\eta(Y)\phi X-g(X,\phi Y)\xi,
\end{equation}
for any $X$ and $Y$. Hence, we also have 
\begin{equation}
\nabla_X\xi=X-\eta(X)\xi.
\end{equation}
Finally, an almost contact metric manifold is said {\it cosymplectic manifold} if and only if $d\eta=0$ and $d\Omega=0$, or equivalently 
\begin{equation}
\nabla\phi=0,
\end{equation}
and in this case, we have
\begin{equation}
\nabla\xi=0.
\end{equation}
The $\phi$-sectional curvature of  an almost contact metric manifold is  defined as the sectional \linebreak curvature on the $2$-planes $\{X,\phi X\}$. When the $\phi$-sectional curvature is constant, we say that the manifold is a space form (Sasakian, Kenmotsu or cosymplectic in each of the three cases above). It is well known that the $\phi$-sectional curvature determines entirely the curvature of the manifold. When the $\phi$-sectional curvature is constant, the curvature tensor is expressed explicitely. Let $R_1^{\star}$, $R_2^{\star}$ and $R_3^{\star}$ be the generalized curvature tensors defined by
\begin{equation}
R_1^{\star}(X,Y)Z=\widetilde{g}(Y,Z)X-\widetilde{g}(X,Z)Y,
\end{equation}
\begin{equation}
R_2^{\star}(X,Y)Z=\eta(X)\eta(Z)Y-\eta(Y)\eta(Z)X+\widetilde{g}(X,Z)\eta(Y)\xi-\widetilde{g}(Y,Z)\eta(X)\xi 
\end{equation} and 
\begin{equation}
R_3^{\star}(X,Y)Z=\Omega(Z,Y)\phi X-\Omega(Z,X)\phi Y+2\Omega(X,Y)\phi Z. 
\end{equation}
For the three cases we are interested in, the curvature of a space form of constant $\phi$-sectional curvature $c$ is given by

\begin{itemize}
\item Sasaki: $R^{\star}=\frac{c+3}{4}R_1^{\star}+\frac{c-1}{4}R_2^{\star}+\frac{c-1}{4}R_3^{\star}.$\\
\item Kenmotsu: $R^{\star}=\frac{c-3}{4}R_1^{\star}+\frac{c+1}{4}R_2^{\star}+\frac{c+1}{4}R_3^{\star}.$\\
\item Cosymplectic: $R^{\star}=\frac{c}{4}R_1^{\star}+\frac{c}{4}R_2^{\star}+\frac{c}{4}R_3^{\star}.$
\end{itemize}
In the sequel, for more clarity, we will denote the Sasakian ({\it resp.}  Kenmotsu, cosymplectic) space form of constant $\phi$-sectional curvature $c$ by $\widetilde{M}_S(c)$ ({\it resp. } $\widetilde{M}_K(c)$, $\widetilde{M}_C(c)$). These space forms appear as particular cases of the so-called generalized Sasakian space forms, introduced by Alegre, Blair and Carriazo in \cite{ABC}. A generalized Sasakian space form, denoted by $\widetilde{M}(f_1,f_2,f_3)$, is a contact metric manifold with curvature tensor of the form
\begin{equation}\label{CurvatureGSasakian}
f_1R_1^{\star}+f_2R_2^{\star}+f_3R_3^{\star},
\end{equation}
where $f_1$, $f_2$ and $f_3$ are real functions on the manifold. The most simple examples of generalized Sasakian space forms are the  warped products of the real line by a complex space form or a generalized complex space forms.  Their conformal deformations as well as their so-called $\mathcal{D}$-homothetic deformations are also generalized Sasakian space forms (see \cite{ABC}). Other examples can be found in \cite{AC}.\\ \\
Now, let $(M,g)$ be a submanifold of an almost contact metric manifold $(\widetilde{M},\widetilde{g},\xi,\eta,\phi)$. The field of tensors $\phi$ induces on $M$, the existence of the following four operators:
$$P:TM\lgra TM,\ N:TM\lgra NM,\ t:NM\lgra TM\ \text{and}\ s:NM\lgra NM ,$$
defined for any $X\in TM$ and $\nu\in NM$. Now, we have
\begin{eqnarray}\label{relationfhst}
\phi X=PX+NX\quad\text{and}\quad
\phi\nu=t\nu+s\nu,
\end{eqnarray} 
where $PX$ and $NX$ are tangential and normal components of $\phi X$, respectively, whereas $t\nu$ and $s\nu$ are the tangential and normal components of $\phi\nu$, respectively. A submanifold $M$ is said invariant ({\it resp.} anti-invariant) if $N$ ({\it resp.} $P$) vanishes identically. In \cite{Lot}, Lotta shows that if the vector field $\xi$ is normal to $M$, then $M$ is anti-invariant.  

\section{\textbf{Biharmonic submanifolds of generalized complex space forms}}\label{sec3}
First of all, we give the following theorem which is a characterization of biharmonic submanifolds in generalized complex space forms.
\begin{thm}\label{thm1}
Let $N(\alpha,\beta)$ be a generalized complex space form and  $M^n$, $n<4$, a submanifold of $N(\alpha,\beta)$ with second fundamental form $B$, shape operator $A$ and mean curvature $H$. Then $M$ is biharmonic if and only if the following two equations are satisfied
$$\left\{
\begin{array}{l}
-\Delta^{\perp}H+\tr \left(B(\cdot,A_H\cdot)\right)-n\alpha H+3\beta klH=0,\\ \\
\frac{n}{2}{\rm grad}|H|^2+2\tr \left(A_{\nabla^{\perp}H}(\cdot)\right)+6\beta jlH=0.
\end{array}
\right.$$
\end{thm}
\noindent
{\bf Proof:} The equations of biharmonicity are well known (see \cite{BMO,Chen,FOR}, for instance). After \linebreak projection of the equation $\tau_2(\psi)=0$ on both tangent and normal bundles, we get the two following equations
\begin{equation}\label{eqbiharmonic}
\left\{
\begin{array}{l}
-\Delta^{\perp}H+\tr \left(B(\cdot,A_H\cdot)\right)+\tr\left(R(\cdot,H)\cdot\right)^{\perp}=0,\\ \\

\frac{n}{2}{\rm grad}|H|^2+2\tr \left(A_{\nabla^{\perp}H}(\cdot)\right)+2\tr\left(R(\cdot,H)\cdot\right)^{\top}=0.
\end{array}
\right.
\end{equation}

Recall that the curvature tensor of $N(\alpha,\beta)$ is given by
\begin{eqnarray*}\label{curvbiharm}
\tr\left(R(\cdot,H)\cdot\right)&=&\alpha\tr\left(R_1(\cdot,H)\cdot\right)+\beta\tr\left(R_2(\cdot,H)\cdot\right).
\end{eqnarray*}
Let us compute the two terms of the right hand side. For this, suppose $\{X_i\}_{i=1}^n$ be a local orthonormal frame of $TM$. First, we have
\begin{eqnarray*}\label{curvbiharm}
\tr\left(R_1(\cdot,H)\cdot\right)&=&\sum_{i=1}^nR_1(X_i,H)X_i.\\
&=&\sum_{i=1}^n\left[g(H,X_i)X_i-g(X_i,X_i)H\right]\\
&=&-nH.
\end{eqnarray*}
Secondly, we have
\begin{eqnarray*}\label{curvbiharm}
\tr\left(R_2(\cdot,H)\cdot\right)&=&\sum_{i=1}^nR_2(X_i,H)X_i.\\
&=&\sum_{i=1}^n\left[g(JH,X_i)JX_i-g(JX_i,X_i)JH+2g(JH,X_i)JX_i\right]\\
&=&3\sum_{i=1}^ng(lH,X_i)JX_i\\
&=&3JlH\\
&=&3jlH+3klH.
\end{eqnarray*}
Since $jlH$ is tangent and $klH$ is normal, by identification of tangent and normal parts, finally, we get the equations of the theorem. \hfill $\square$\\

We have the following corollary for the particular cases of hypersurfaces, Lagrangian or complex surfaces and curves.
\begin{cor}\label{cor1}
Let $N(\alpha,\beta)$ be a generalized complex space form and $M$ a submanifold of $N(\alpha,\beta)$ with second fundamental form $B$, shape operator $A$ and mean curvature $H$.
\begin{enumerate}
\item If $M$ is a hypersurface, then $M$ is biharmonic if and only if
$$\left\{
\begin{array}{l}
-\Delta^{\perp}H+\tr \left(B(\cdot,A_H\cdot)\right)-3(\alpha+\beta) H=0,\\ \\
\frac{3}{2}{\rm grad}|H|^2+2\tr \left(A_{\nabla^{\perp}H}(\cdot)\right)=0.
\end{array}
\right.$$
\item If $M$ is a complex surface, then $M$ is biharmonic if and only if
$$\left\{
\begin{array}{l}
-\Delta^{\perp}H+\tr \left(B(\cdot,A_H\cdot)\right)-2\alpha H=0,\\ \\
{\rm grad}|H|^2+2\tr \left(A_{\nabla^{\perp}H}(\cdot)\right)=0.
\end{array}
\right.$$
\item If $M$ is a Lagrangian surface, then $M$ is biharmonic if and only if
$$\left\{
\begin{array}{l}
-\Delta^{\perp}H+\tr \left(B(\cdot,A_H\cdot)\right)-2\alpha H-3\beta H=0,\\ \\
{\rm grad}|H|^2+2\tr \left(A_{\nabla^{\perp}H}(\cdot)\right)=0.
\end{array}
\right.$$
\item If $M$ is a curve, then $M$ is biharmonic if and only if
$$\left\{
\begin{array}{l}
-\Delta^{\perp}H+\tr\left( B(\cdot,A_H\cdot)\right)-\alpha H-3\beta (H+m^2H)=0,\\ \\
\frac{1}{2}{\rm grad}|H|^2+2\tr \left(A_{\nabla^{\perp}H}(\cdot)\right)=0.
\end{array}
\right.$$
\end{enumerate}
\end{cor}
{\bf Proof:} The proof is a direct consequence of Theorem \ref{thm1}.
\begin{enumerate}
\item If $M$ is a hypersurface, then $J$ maps normal vectors on tangent vectors, that is, $m=0$. Hence, by relation \eqref{relation1.2}, we have $klH=-H$ and by relation \eqref{relation1.3}, $jlH=0$, which gives the result by Theorem \ref{thm1}.
\\
\item If $M$ is a complex surface, then, $k=0$ and $l=0$.\\
\item If $M$ is a Lagrangian surface, then $j=0$ and $m=0$. Moreover, since $m=0$, as for hypersurfaces, we have $klH=-H$ by relation \eqref{relation1.2}.\\
\item If $M$ is a curve, then $j=0$. Hence, by relation \eqref{relation1.2}, $klH=-(H+m^2H)$.

\end{enumerate}
\hfill$\square$

\begin{remark}
It is a well known fact that any complex submanifold of a K\"ahler manifold is \linebreak necessarily minimal. But as mentioned above, the generalized space forms $N(\alpha,\beta)$ are not K\"ahlerian unless there are the complex projective plane or the complex hyperbolic plane. Hence, considering biharmonic surfaces into $N(\alpha,\beta)$ is of real interest, since they are not necessarily minimal.
\end{remark}

We obtain immediately the following corollaries for curves and complex or Lagranian surfaces with parallel mean curvature. 
\begin{cor}\label{corlag}
\begin{enumerate}
\item If $M$ be a Lagrangian surface of $N(\alpha,\beta)$ with parallel mean curvature, then $M$ is biharmonic if and only if
$$\tr\left(B(\cdot,A_H\cdot)\right)=(2\alpha+3\beta)H.$$
\item If $M$ be a complex surface of $N(\alpha,\beta)$ with parallel mean curvature, then $M$ is biharmonic if and only if
$$\tr\left(B(\cdot,A_H\cdot)\right)=2\alpha H.$$
\item If $M$ is a curve in $N(\alpha,\beta)$ with parallel mean curvature, then $M$ is biharmonic if and only if
$$\tr\left( B(\cdot,A_H\cdot)\right)=\alpha H+3\beta (H+m^2H)=0.$$
\end{enumerate}
\end{cor}
Now, we give some curvature properties of constant mean curvature submanifolds in $N(\alpha,\beta)$. We have this first proposition for hypersurfaces.
\begin{prop}\label{propB}
Let $N(\alpha,\beta)$ be a generalized complex space form and $M^3$ a hypersurface of $N(\alpha,\beta)$ with non zero constant mean curvature $H$. Then, $M$ is proper-biharmonic if and only if
 $$||B||^2=3(\alpha+\beta),$$
or equivalentely, if the scalar curvature of $M$ satisfies
$$\scal_M=3(\alpha+\beta)+9H^2.$$

\end{prop}
\begin{remark}
In particular, the norm of the second fundamental form and the scalar curvature of $M$ are constant. \end{remark}
\noindent
{\bf Proof:} As $M$ is a hypersurface, by Corollary \ref{cor1}, $M$ is biharmonic if and only if
$$\left\{
\begin{array}{l}
-\Delta^{\perp}H+\tr B(\cdot,A_H\cdot)-3\alpha H-3\beta H=0,\\ \\
\frac{3}{2}{\rm grad}|H|^2+2\tr A_{\nabla^{\perp}H}(\cdot)=0.
\end{array}
\right.$$
Since $M$ has constant mean curvature, the second equation is trivially satisfied and the first becomes
$$\tr \Big(B(\cdot,A_H(\cdot))\Big)=(3\alpha+3\beta)H.$$
Moreover, for hypersurfaces, we have $A_H=HA$ which implies
$$\tr \Big(B(\cdot,A_H(\cdot))\Big)=H\tr \Big(B(\cdot,A(\cdot))\Big)=H||B||^2.$$
Finally, since $H$ is a non-zero constant, we get the desired identity $|B|^2=3(\alpha+\beta)$.\\
For the second equivalence, by the Gauss equation, we have
$$
\scal_M=\sum_{i,j=1}^3g\left( R^N(X_i,X_j)X_j,X_i\right)-||B||^2+9H^2,
$$
where $\{X_1,X_2,X_3\}$ is a local orthonormal frame of $M$. From the expression of the curvature tensor of $N(\alpha,\beta)$, we get
$$
\scal_M=6(\alpha+\beta)-||B||^2+9H^2.
$$
Hence, we deduce that $M$ is proper biharmonic if and only if $||B||^2=3(\alpha+\beta)$, that is, if and only if $\scal_M=3(\alpha+\beta)+9H^2$.
\hfill$\square$\\ \\
An immediate consequence of this proposition is the following corollary.
\begin{cor}
There exists no biharmonic hypersurface with constant mean curvature in a generalized complex space form $N(\alpha,\beta)$ of negative scalar curvature.
\end{cor}
{\bf Proof:} From Proposition \ref{propB}, a constant mean curvature hypersurface of $N(\alpha,\beta)$ is biharmonic if and only if $|B|^2=3(\alpha+\beta)$, which is possible only if $\alpha+\beta$ is positive, that is if $N(\alpha,\beta)$ has positive scalar curvature.
\hfill $\square$\\ \\
Finally, we give this last proposition which give an estimate of the mean curvature for a biharmonic Lagrangian surface. 
\begin{prop}\label{proplag}
\begin{enumerate} 
\item There exists no proper biharmonic Lagrangian surface with constant mean curvature in $N(\alpha,\beta)$ if the $2\alpha+3\beta$ is non-positive everywhere. \\
\item
Suppose that $2\alpha+3\beta$ is a positive function. Let $M$ be a Lagrangian surface of  $N(\alpha,\beta)$ with non-zero constant mean curvature. Then we have the following observations.
\begin{enumerate} 
\item If $M$ is proper-biharmonic, then $0<|H|^2 \leqslant \inf_{M}\left(\frac{2\alpha+3\beta}{2}\right)$.  
\item If $|H|^2= \inf_M\left(\frac{2\alpha+3\beta}{2}\right)$, then $M$ is proper-bihramonic if and only if $\alpha$ and $\beta$ are constant over $M$, $M$ is pseudo-umbilical and $\nabla^{\perp}H=0$.
\end{enumerate}
\end{enumerate}
\end{prop}
{\bf Proof:} Let $M$ be a biharmonic submanifold of $N(\alpha,\beta)$ with non-zero constant mean curvature. Since $M$ is a Lagrangian surface, by the third assertion of Corollary \ref{cor1}, we have
$$-\Delta^{\perp}H+\tr \left(B(\cdot,A_H\cdot)\right)-(2\alpha+3\beta)H=0.$$
Hence, by taking the scalar product with $H$, we have
$$-\left\langle\Delta^{\perp}H,H\right\rangle=(2\alpha+3\beta)|H|^2-\tr B(\cdot,A_H\cdot).$$
Using the Bochner formula and the fact that $|H|$ is constant, we get
$$(2\alpha+3\beta)|H|^2=\tr \left(B(\cdot,A_H\cdot)\right)+|\nabla^{\perp}H|^2.$$
Moreover, by Cauchy-Schwarz inequality, we get $\tr \left(B(\cdot,A_H\cdot)\right)\geqslant 2|H|^4$. Therefore, we have 
$$(2\alpha+3\beta)|H|^2\geqslant 2|H|^4+|\nabla^{\perp}H|^2\geqslant 2|H|^4.$$
Since $|H|$ is a non-zero constant, we have 
$0<H^2\leqslant \inf_M\left(\frac{2\alpha+3\beta}{2}\right)$. This is only possible if the function $2\alpha+3\beta$ is positive everywhere. This remark gives the non-existence result (1) and the first point of (2).\\ \\
Now, assume that $|H|^2= \inf_M\left(\frac{2\alpha+3\beta}{2}\right)$. If $M$ is proper biharmonic, then all the inequalites\linebreak above become equality. First, $2\alpha+3\beta$ is a constant. But $\alpha+\beta$ is also a constant, then $\alpha$ and $\beta$ are constant. Hence, necessarily, $\alpha=\beta$ over $M$. Moreover, equality occurs in the Cauchy-Schwarz inequality, i.e., $M$ is pseudo-umbilical. Finally, we also have $\nabla^{\perp}H=0$. \\
Conversly, if $\alpha$ and $\beta$ are constant over $M$, $M$ is pseudo-umbilical and $\nabla^{\perp}H=0$, then we have $2\alpha+3\beta=|H|^2$ and we get immediatley
$$-\Delta^{\perp}H+\tr \left(B(\cdot,A_H\cdot)\right)-(2\alpha+3\beta)H=0.$$
and
$${\rm grad}|H|^2+2\tr A_{\nabla^{\perp}H}(\cdot)=0.$$
Hence, by Corollary \ref{cor1}, $M$ is biharmonic. This concludes the proof.
\hfill $\square$\\ \\
We have an analogous result for complex surfaces. Note again that, in this context, complex\linebreak surfaces are not necessarily minimal.
\begin{prop}\label{propcomp}
\begin{enumerate} 
\item There exists no proper biharmonic complex surface with constant mean curvature in $N(\alpha,\beta)$ if the function $\alpha$ is non-positive everywhere. \\
\item
Suppose that $\alpha$ is a positive function.
Let $M$ be a complex surface of  $N(\alpha,\beta)$ ($\alpha+\beta>0$) with non-zero constant mean curvature. Then we have
\begin{enumerate} 
\item If $M$ is proper-biharmonic, then $0<|H|^2 \leqslant \inf_{M}\left(\alpha\right)$.  
\item If $|H|^2= \inf_M(\alpha)$, then $M$ is proper-bihramonic if and only if $\alpha$ and $\beta$ are constant over $M$, $M$ is pseudo-umbilical and $\nabla^{\perp}H=0$.
\end{enumerate}
\end{enumerate}
\end{prop}
{\bf Proof:} Let $M$ be a biharmonic submanifold of $N(\alpha,\beta)$ with non-zero constant mean curvature. Since $M$ is a complex surface, by the second assertion of Corollary \ref{cor1}, we have
$$-\Delta^{\perp}H+\tr \left(B(\cdot,A_H\cdot)\right)-2\alpha H=0.$$
The rest of the proof is analogous to the proof in the Lagrangian case with $2\alpha$ instead of $2\alpha+3\beta$.
\hfill $\square$\\ \\
Note that the results of this section contains the particular case of the complex projective\linebreak planes $\CC P^2(4\alpha)$, proved in \cite{FLMO}.

\section{\textbf{Biharmonic submanifolds of generalized Sasakian space forms}}\label{sec4}
Now, we consider biharmonic submanifolds of generalized Sasakian space forms. First we give the necessary and sufficient condition for a submanifold of a generalized Sasakian space form to be biharmonic.
\begin{thm}\label{thm2}
Let $(M^n,g)$ be a submanifold of a generalized Sasakian space form $\widetilde{M}(f_1,f_2,f_3)$ with second fundamental form $B$, shape operator $A$ and mean curvature $H$. Then $M$ is biharmonic if and only if both equations are fulfilled:
$$\left\{
\begin{array}{l}
-\Delta^{\perp}H+\tr B(\cdot,A_H\cdot)=nf_1H-f_2|\xi^{\top}|^2H-nf_2\eta(H)\xi^{\perp}-3f_3NtH,\\ \\
\frac{n}{2}\grad|H|^2+2\tr A_{\nabla^{\perp}H}(\cdot)=-2f_2(n-1)\eta(H)\xi^{\top}-6f_3PtH.
\end{array}
\right.$$
\end{thm}
{\bf Proof:} The curvature tensor of generalized Sasakian space form is given by equation \ref{CurvatureGSasakian}. Now, we have

\begin{eqnarray*}
R^{\star}(X, Y)Z &=& f_{1}R_{1}^{\star}(X, Y)Z + f_{2}R_{2}^{\star}(X, Y)Z + f_{3}R_{3}^{\star}(X, Y)Z \\
&=& f_{1}\{\tilde{g}(Y, Z)X - \tilde{g}(X, Z)Y\}\\
&+& f_{2}\{\eta(X)\eta(Z)Y - \eta(Y)\eta(Z)X + \tilde{g}(X, Z)\eta(Y)\xi - \tilde{g}(Y, Z)\eta(X)\xi \} \\ &+& f_{3}\{\tilde{g}(X, \phi Z)\phi Y - \tilde{g}(Y, \phi Z)\phi X + 2\tilde{g}(X, \phi Y)\phi Z\}.
\end{eqnarray*}

From this equation, we have

\begin{eqnarray*}
R^{\star}(X_{i}, H)X_{i} &=& f_{1}\{\tilde{g}(H, X_{i})X_{i} - \tilde{g}(X_{i}, X_{i})H\}
+ f_{2}\{\eta(X_{i})\eta(X_{i})H - \eta(H)\eta(X_{i})X_{i} + \tilde{g}(X_{i}, X_{i})\eta(H)\xi \} \\ &+& f_{3}\{\tilde{g}(X_{i}, \phi X_{i})\phi H - \tilde{g}(H, \phi X_{i})\phi X_{i} + 2\tilde{g}(X_{i}, \phi H)\phi X_{i}\}.
\end{eqnarray*}

From \eqref{relationfhst}, we get 

\begin{eqnarray*}
tr R^{\star}(X_{i}, H)X_{i} &=&  - f_{1}nH 
+ f_{2}\sum_{i}\{\eta(X_{i})^{2}H - \eta(H)\eta(X_{i})X_{i} + |X_{i}|^{2}\eta(H)\xi \} \\ &+& f_{3}\sum_{i}\{tr(P)\phi H - \tilde{g}(H, N X_{i})\phi X_{i} + 2\tilde{g}(X_{i}, t H)\phi X_{i}\}
\\
&=&- f_{1}nH 
+ f_{2}\{|\xi^{\top}|^{2}H - \eta(H)\xi^{\top} + n\eta(H)\xi\} \\ &+& f_{3}\sum_{i}\{tr(P)t H + tr(P)s H - \tilde{g}(H, N X_{i})P X_{i} - \tilde{g}(H, N X_{i})N X_{i}\\ &+& 2\tilde{g}(X_{i}, t H)P X_{i} + 2\tilde{g}(X_{i}, t H)N X_{i}\}.
\end{eqnarray*}

By the anti-symmetry of $\phi$, we have immediately $tr(P)=0$ and using the fact $$\tilde{g}(H, N X_{i}) = - \tilde{g}(t H, X_{i}),$$ we get
\begin{eqnarray*}
tr R^{\star}(X_{i}, H)X_{i} &=& - f_{1}nH 
+ f_{2}\{|\xi^{\top}|^{2}H - \eta(H)\xi^{\top} + n\eta(H)\xi\}+3f_3(PtH+NtH).
\end{eqnarray*}
Finally, reporting in Equation \eqref{eqbiharmonic} the tangential and normal parts, we have result of the theorem.
\hfill$\square$\\\\
From this general condition, we can state many interesting particular cases. Namely, we have the following corollary.
\begin{cor}\label{cor2}
Let $(M^n,g)$ be a submanifold of a generalized Sasakian space form $\widetilde{M}(f_1,f_2,f_3)$ with second fundamental form $B$, shape operator $A$ and mean curvature $H$.
\begin{enumerate}
\item If $M$ is invariant, then $M$ is biharmonic if and only if
$$\left\{
\begin{array}{l}
-\Delta^{\perp}H+\tr B(\cdot,A_H\cdot)=nf_1H-f_2|\xi^{\top}|^2H-nf_2\eta(H)\xi^{\perp},\\ \\
\frac{n}{2}\grad|H|^2+2\tr A_{\nabla^{\perp}H}(\cdot)=-2f_2(n-1)\eta(H)\xi^{\top}-6f_3PtH.
\end{array}
\right.$$
\item  If $M$ is anti-invariant, then $M$ is biharmonic if and only if
$$\left\{
\begin{array}{l}
-\Delta^{\perp}H+\tr B(\cdot,A_H\cdot)=nf_1H-f_2|\xi^{\top}|^2H-nf_2\eta(H)\xi^{\perp}-3f_3NtH,\\ \\
\frac{n}{2}\grad|H|^2+2\tr A_{\nabla^{\perp}H}(\cdot)=-2f_2(n-1)\eta(H)\xi^{\top}.
\end{array}
\right.$$
\item
If $\xi$ is normal to $M$, then $M$ is biharmonic if and only if
$$\left\{
\begin{array}{l}
-\Delta^{\perp}H+\tr B(\cdot,A_H\cdot)=nf_1H-nf_2\eta(H)\xi-3f_3NtH.\\ \\
\frac{n}{2}\grad|H|^2+2\tr A_{\nabla^{\perp}H}(\cdot)=0
\end{array}
\right.$$
\item If $\xi$ is tangent to $M$, then $M$ is biharmonic if and only if
$$\left\{
\begin{array}{l}
-\Delta^{\perp}H+\tr B(\cdot,A_H\cdot)=nf_1H-f_2H-3f_3NtH,\\ \\
\frac{n}{2}\grad|H|^2+2\tr A_{\nabla^{\perp}H}(\cdot)=-6f_3PtH.
\end{array}
\right.$$
\item If $M$ is a hypersurface, then $M$ is biharmonic if and only if
$$\left\{
\begin{array}{l}
-\Delta^{\perp}H+\tr B(\cdot,A_H\cdot)=(nf_1+3f_3)H-f_2|\xi^{\top}|^2H-(nf_2+3f_3)\eta(H)\xi^{\perp}\\ \\
\frac{n}{2}\grad|H|^2+2\tr A_{\nabla^{\perp}H}(\cdot)=-(2(n-1)f_1+6f_3)\eta(H)\xi^{\top}.
\end{array}
\right.$$
\end{enumerate}
\end{cor}
{\it \textbf{Proof:}} The proof is a direct consequence of the above theorem, using the following facts:
\begin{enumerate}
\item $P=0$ for $M$ invariant.\\
\item $N=0$ for $M$ anti-invariant.\\
\item $\xi^{\top}=0$ and $\xi^{\perp}=\xi$ if $\xi$ is normal. Moreover, since $\xi$ is normal, $M$ is necessarily anti-invariant and so $P=0$.\\
\item $\xi^{\top}=\xi$, $|\xi|=1$ and $\xi^{\perp}=0$ if $\xi$ is tangent.\\
\item If $M$ is a hypersurface, then $\phi(H)$ is tangent and so $sH=0$. From this, we get
$$-H+\eta(H)\xi=\phi^2H=PtH+NtH.$$
By identification of tangential and normal parts, we obtain 
$$NtH=-H+\eta(H)\xi^{\perp}\quad\text{and}\quad PtH=\eta(H)\xi^{\top},$$
which gives the result.\\
\end{enumerate}
Now, using these characterizations of biharmonic submanifolds, we can obtain some curvature properties in some special case. First, analogously to the case of generalized complex space forms (Proposition \ref{propB}), we have the following result for hypersurfaces.
\begin{prop}\label{propscal}
Let $(M^n,g)$ by a hypersurface of a generalized Sasakian space form $\widetilde{M}(f_1,f_2,f_3)$ of non zero constant mean curvature and suppose that $\xi$ is tangent to $M$. Then $M$ is proper biharmonic if and only if
$$|B|^2=nf_1-f_2+3f_3$$
or equivalently if and only if
$$\scal_M=n(n-2)f_1+(2n-1)f_2-3nf_3+(n-1)H^2.$$
\end{prop}
Before proving this proposition, we give the following lemma.
\begin{lemma}\label{xitangent}
If $M$ is a hypersurface of an almost contact metric manifold with $\xi$ tangent to $M$, then $Pt=0$ and $Nt=-{\rm Id}$.
\end{lemma}
{\bf Proof:} Since $\xi$ is tangent to $M$, we have $\eta(\nu)=\langle\xi,\nu\rangle=0$ and so
$$\phi^2\nu=-\nu+\eta(\nu)\xi=-\nu.$$
On the other hand, we have
\begin{eqnarray*}
\phi^2\nu&=&\phi(t\nu+s\nu)\\
&=&Pt\nu+Nt\nu+ts\nu+s^2\nu
\end{eqnarray*}
Hence, we get
\begin{equation}\label{minusnu}
-\nu=Pt\nu+Nt\nu+ts\nu+s^2\nu.
\end{equation}
Moreover, since $\langle\phi\nu,\nu\rangle=\Omega(\nu,\nu)=0$, we have that $\phi\nu$ is tangent, i.e., $s\nu=0$. Thus, Equation \eqref{minusnu} becomes
$$-\nu=Pt\nu+Nt\nu,$$
and so $Pt=0$ and $Nt=-{\rm Id}$ by identification of tangential and normal parts.
\hfill$\square$

{\bf Proof of Proposition \ref{propscal}:} Since $M$ is a hypersurface of $\widetilde{M}(f_1,f_2,f_3)$ with non zero constant mean curvature and with $\xi$ tangent to $M$, we know by Corollary \ref{cor2} that $M$ is biharmonic if and only if
$$\left\{
\begin{array}{l}
\tr B(\cdot,A_H\cdot)=nf_1H-f_2H-3f_3NtH,\\ \\
2f_2(n-1)\eta(H)\xi+6f_3PtH=0.
\end{array}
\right.$$
Moreover, $\eta(H)=\langle H,\xi\rangle=0$ since $H$ is normal and $\xi$ tangent. In addition, by Lemma \ref{xitangent}, we have $Pt=0$ and $Nt=-{\rm Id}$, hence, the second equation is trivial and the first becomes
$$\tr B(\cdot,A_H\cdot)=nf_1H-f_2H+3f_3H,$$
or equivalently
$$|B|^2=nf_1-f_2+3f_3,$$
since $\tr B(\cdot,A_H\cdot)=|B|^2H$ and $H$ is a non zero constant.\\
\\
Now for second part, we have by the Gauss formula
\begin{eqnarray*}
Scal_{M} &=& \sum_{i, j}\tilde{g}(R^{\star}(X_{i}, X_{j})X_{j}, X_{i}) - |B|^{2} - nH^{2}\\
&=& \sum_{i, j}f_{1}\{\tilde{g}(X_{j}, X_{j})\tilde{g}(X_{i}, X_{i}) - \tilde{g}(X_{i}, X_{j})\tilde{g}(X_{j}, X_{i})\}
+ \sum_{i, j}f_{2}\{\eta(X_{i})\eta(X_{j})\tilde{g}(X_{j}, X_{i})\\ &-& \eta(X_{j})\eta(X_{j})\tilde{g}(X_{i}, X_{i}) + \tilde{g}(X_{i}, X_{j})\eta(X_{j})\tilde{g}(\xi, X_{i}) - \tilde{g}(X_{j}, X_{j})\eta(X_{i})\tilde{g}(\xi, X_{i})\}\\
&+& \sum_{i, j}f_{3}\{\tilde{g}(X_{i}, \phi X_{j})\tilde{g}(\phi X_{j}, X_{i}) - \tilde{g}(X_{j}, \phi X_{j})\tilde{g}(\phi X_{i}, X_{i}) + 2\tilde{g}(X_{i}, \phi X_{j})\tilde{g}(\phi X_{j}, X_{i})\}\\ &-& |B|^{2} - nH^{2}\\
&=& n(n - 1)f_{1} + 2(n - 1)f_{2} - (n - 1)f_{3} - |B|^{2} - nH^{2}.
\end{eqnarray*}
Using the value of $|B|^{2}$ obtain in the first part of the proof, we get the required result, that is,
$$\scal_M=n(n-2)f_1+(2n-1)f_2-3nf_3+(n-1)H^2.$$
This concludes the proof.
\hfill$\square$\\ \\
Now, from this proposition, we can prove the following non-existence result.
\begin{cor}
There exists no proper biharmonic CMC hypersurface with $\xi$ tangent in a \linebreak generalized Sasakian space form $\widetilde{M}(f_1,f_2,f_3)$ (of dimension $n+1$) if the functions $f_1, f_2, f_3$ satisfy $nf_1-f_2+3f_3\leqslant 0$. In particular, there exists no proper biharmonic CMC hypersurface with $\xi$ tangent in
\begin{itemize}
\item a Sasaki space form $\widetilde{M}^{n+1}_S(c)$\hspace{.2cm} if \hspace{.2cm} $c\leqslant-\frac{3n-2}{n+2}$,
\item a Kenmotsu space form $\widetilde{M}^{n+1}_K(c)$\hspace{.2cm} if\hspace{.2cm} $c\leqslant\frac{3n-2}{n+2}$,
\item a cosymplectic space form $\widetilde{M}^{n+1}_C(c)$\hspace{.2cm} if\hspace{.2cm} $c\leqslant0$.
\end{itemize}
\end{cor}
{\bf Proof:}
We know from Proposition \ref{propscal} that a hypersurface $M$ of $\widetilde{M}(f_1,f_2,f_3)$ with non zero \linebreak constant mean curvature $H$ and $\xi$ tangent to $M$ is biharmonic if and only if its second \linebreak fundamental form $B$ satisfies
$$|B|^2=nf_1-f_2+3f_3.$$\\
In particular, this is not possible if 
\begin{equation}\label{condf1f2f3}
nf_1-f_2+3f_3\leqslant0.
\end{equation}\\
If $\widetilde{M}(f_1,f_2,f_3)$ is a Sasakian space form of $\phi$-sectional curvature $c$, we have $f_1=\frac{c+3}{4}$ and $f_2=f_3=\frac{c-1}{4}$. Therefore, the assumption $nf_1-f_2+3f_3\leqslant0$ reduces to $c\leqslant-\frac{3n-2}{n+2}$. For the Kenmotsu ({\it resp. } cosymplectic) case, we have $f_1=\frac{c-3}{4}$ and $f_2=f_3=\frac{c+1}{4}$ ({\it resp.} $f_1=f_2=f_3=\frac c4$) and the condition $nf_1-f_2+3f_3\leqslant0$ becomes $c\leqslant\frac{3n-2}{n+2}$ ({\it resp.} $c\leqslant0$).
\hfill $\square$\\ \\
For the next two results, we introduce for more clarity the following constant. 
$$ K(m,c)=\left\{\begin{array}{ll}
\frac{(m+2)c}{4}+\frac{(3m-2)}{4}&\text{if}\ M(f_1,f_2,f_3)= \widetilde{M}^{n+1}_S(c),\\ \\
\frac{(m+2)c}{4}-\frac{(3m-2)}{4}&\text{if}\ M(f_1,f_2,f_3)= \widetilde{M}^{n+1}_K(c),\\ \\
\frac{(m+2)c}{4}&\text{if}\ M(f_1,f_2,f_3)= \widetilde{M}^{n+1}_C(c).
\end{array}
\right.$$
We remark that this constant is just the value of $mf_1-f_2+3f_3$ for the corresponding functions $f_1$, $f_2$ and $f_3$ of each space form $\widetilde{M}^{n+1}_S(c)$, $\widetilde{M}^{n+1}_K(c)$ and $\widetilde{M}^{n+1}_C(c)$. \\ \\
We are now able to state the following propositions which are the counterpart in this almost\linebreak contact setting of Propositions \ref{proplag} and \ref{propcomp}. Here Lagrangian surface and complex surface are replaced, respectively by $\xi$ is normal and $\xi$ is tangent. Therefore, we have:
\begin{prop}
\begin{enumerate} 
\item There exists no proper biharmonic  of constant mean curvature $H$ so that $\xi$ and $\phi H$ are tangent in  $\widetilde{M}^{n+1}_S(c)$ ({\it resp.} $\widetilde{M}^{n+1}_K(c)$ or $\widetilde{M}^{n+1}_C(c)$) of constant mean curvature $H$ so that $\xi$ and $\phi H$ are tangent if $K(m,c)\leqslant 0$..\\
\item
Assume that $K(m,c)>0$. Let $(M^m,g)$ a submanifold of  $\widetilde{M}^{n+1}_S(c)$ ({\it resp.} $\widetilde{M}^{n+1}_K(c)$ or $\widetilde{M}^{n+1}_C(c)$) of constant mean curvature $H$ so that $\xi$ and $\phi H$ are tangent. Then
\begin{enumerate}
\item If $M$ is proper biharmonic, then $|H|^2\in \left(0,\frac{K(m,c)}{m}\right]$.
\item If $|H|^2=\frac{K(m,c)}{m}$, then $M$ is proper biharmonic if and only if $M$ is pseudo-umbilical and $\nabla^{\perp}H=0$.
\end{enumerate}
\end{enumerate}
\end{prop}
{\bf Proof:} Since $\phi H$ is tangent, we have $sH=0$. Hence, we deduce that $\phi^2H=PtH+NtH$. But, we also have $\phi^2H=-H+\eta(H)\xi$. Now, since $\xi$ is tangent, we have $\eta(H)=0$ and so $\phi^2H=-H$. Identifying tangential and normal parts, we get $PtH=0$ and $NtH=-H$. Moreover, $M$ is proper biharmonic, so we get form Corollary \ref{cor2} that \begin{eqnarray*}-\Delta^{\perp}H+\tr B(\cdot,A_H\cdot)&=&mf_1H-f_2H+3f_3H\\
&=&K(m,c)H.
\end{eqnarray*}
Now, we take the scalar product by $H$ to obtain
 $$-\left\langle\Delta^{\perp}H,H\right\rangle+\left\langle\tr B(\cdot,A_H\cdot),H\right\rangle=K(m,c)|H|^2.$$
Form the B\"ochner formula, we get
$$\frac{1}{2}\Delta|H|^2=\left\langle\Delta^{\perp}H,H\right\rangle-|\nabla^{\perp}H|^2$$
Moreover, from the fact that $\left\langle\tr B(\cdot,A_H\cdot),H\right\rangle=|A_H|^2$ and $|H|$ is constant, we get 
$$|A_H|^2+|\nabla^{\perp}H|^2=K(m,c)|H|^2.$$
Finally, by the Cauchy-Schwarz inequality, we have $|A_H|^2\geqslant \frac{1}{m}\tr(A_H)=m|H|^4$, which gives
$$K(m,c)|H|^2=|A_H|^2+|\nabla^{\perp}H|^4\geqslant m|H|^2+|\nabla^{\perp}H|^2\geqslant m|H|^4.$$
Since $|H|$ is a positive constant, we obtain
$$K(m,c)\geqslant m|H|^2,$$
which concludes the proof of the first point (this is not possible is $K(m,c)\leqslant0$ and the first assertion of point (2).\\
For the second part, we assume that $K(m,c)= m|H|^2$. Since $M$ is biharmonic, then all the inequalities in the proof of part $1$ become equalities. In particular, we have $\nabla^{\perp}H=0$ and equality occurs in the Cauchy-Schwarz inequality, that is, $A_H$ is scalar or in other terms, $M$ is pseudo-umbilical.\\
Conversely, if $\nabla^{\perp}H=0$ and $M$ is pseudo umbilical, then, we see immediatley that both equations of Theorem \ref{thm1} are fulfilled and so $M$ is proper biharmonic. Indeed, the first equality is achieved by the above discussion and the second is trivial since $|H|^2$ is constant, $\nabla^{\perp}H=0$ and as we have already seen, $PtH=0$.
\hfill$\square$
\begin{remark}\label{rem1}
\begin{enumerate}
\item The assumption that $\phi H$ is tangent is automatically verified for \linebreak hypersurfaces with non zero constant mean curvature. Indeed, since $\phi$ is anti-symmetric, we have $<\phi H,H>=0$ and thus $\phi H$ is tangent. Of course, this fact is specific to codimension 1.
\item
Note that the condition $K(m,c)>0$ is
$$ \left\{\begin{array}{ll}
c>\frac{2-3m}{m+2}&\text{if}\ M(f_1,f_2,f_3)= \widetilde{M}^{n+1}_S(c),\\ \\
c>\frac{3m-2}{m+2}&\text{if}\ M(f_1,f_2,f_3)= \widetilde{M}^{n+1}_K(c),\\ \\
c>0&\text{if}\ M(f_1,f_2,f_3)= \widetilde{M}^{n+1}_C(c).
\end{array}
\right.$$
\end{enumerate}
\end{remark}

\begin{prop}
\begin{enumerate} 
\item There exists no proper biharmonic submanifolds of constant mean curvature $H$ so that $\xi$ is tangent and $\phi H$ is normal in $\widetilde{M}^{n+1}_S(c)$ ({\it resp.} $\widetilde{M}^{n+1}_K(c)$ or $\widetilde{M}^{n+1}_C(c)$)  if $K(m,c)\leqslant3$. \\
\item
Let $(M^m,g)$ a submanifold of  $\widetilde{M}^{n+1}_S(c)$ ({\it resp.} $\widetilde{M}^{n+1}_K(c)$ or $\widetilde{M}^{n+1}_C(c)$) of constant mean curvature $H$ so that $\xi$ is tangent and $\phi H$ is normal. Suppose that $K(m,c)>3$. Then, we have
\begin{enumerate}
\item If $M$ is proper biharmonic, then $|H|^2\in \left(0,\frac{K(m,c)-3}{m}\right]$.
\item If $|H|^2=\frac{K(m,c)-3}{m}$, then $M$ is proper biharmonic if and only if $M$ is pseudo-umbilical and $\nabla^{\perp}H=0$.
\end{enumerate}
\end{enumerate}
\end{prop}
{\bf Proof:} The proof is comparable to the proof of the previous theorem. Here, since $\phi H$ is normal, we have $tH=0$. Moreover, since $\xi$ is tangent and $M$ is proper biharmonic, we get from Corollary \ref{cor2} that \begin{eqnarray*}-\Delta^{\perp}H+\tr B(\cdot,A_H\cdot)&=&mf_1H-f_2H\\
&=&(K(m,c)-3)H.
\end{eqnarray*}
As for the previous lemma, taking the scalar product by $H$ and using the B\"ochner formula and the Cauchy-Schwarz inequality, we get
$$(K(m,c)-3)|H|^2=|A_H|^2+|\nabla^{\perp}H|^4\geqslant m|H|^2+|\nabla^{\perp}H|^2\geqslant m|H|^4.$$
Since $|H|$ is a positive constant, we obtain $(K(m,c)-3)\geqslant m|H|^2$ which concludes the proof of the first assertion and the first poitn of the second assertion.\\
Now, if $K(m,c)-3=|H|^2$, then $M$ is biharmonic if and only if the above inequalities are equalities, that is $\nabla^{\perp}H=0$ and $M$ is pseudo umbilical.
\hfill$\square$
\begin{remark}
\begin{enumerate}
\item
The assumption that $\phi H$ is normal, can not be verified for hypersurfaces by Remark \ref{rem1}. Hence, this proposition holds only for codimension at least $2$.
\item
It is to note that the condition $K(m,c)>3$ is
$$ \left\{\begin{array}{ll}
c>\frac{14-3m}{m+2}&\text{if}\ M(f_1,f_2, f_3)= \widetilde{M}^{n+1}_S(c),\\ \\
c>\frac{10+3m}{m+2}&\text{if}\ M(f_1,f_2, f_3)= \widetilde{M}^{n+1}_K(c),\\ \\
c>\frac{12}{m+2}&\text{if}\ M(f_1,f_2, f_3)= \widetilde{M}^{n+1}_C(c).
\end{array}
\right.$$
\end{enumerate}
\end{remark}
\section*{\textbf{Acknowledgements}}
Second author is supported by post doctoral scholarship of ``Harish Chandra Research Institute", Department of Atomic Energy, Government of India.

\end{document}